\documentclass[12pt]{amsart}

\usepackage{amssymb,bm,mathrsfs,mathtools,booktabs}

\usepackage{enumerate}
\usepackage{tikz}
\usepackage[centering,width=6in]{geometry}
\usepackage{graphics,verbatim}
\usepackage{graphicx}
\usepackage{csquotes}
\usepackage{hyperref}
\hypersetup{
    colorlinks,
    citecolor=blue,
    filecolor=black,
    linkcolor=red,
    urlcolor=black
}
\usepackage{enumitem}
\usepackage{thmtools} 
\usepackage{thm-restate}
\usepackage{dsfont}

\newtheorem{theorem}{Theorem}[section]
\newtheorem{lemma}[theorem]{Lemma}
\newtheorem{proposition}[theorem]{Proposition}

\newtheorem{conj}[theorem]{Conjecture}

\theoremstyle{definition}
\newtheorem{definition}[theorem]{Definition}
\newtheorem{example}[theorem]{Example}

\theoremstyle{remark}

\newtheorem{ques}{Question}

\numberwithin{equation}{section}

\newcommand{\R}{{\mathbb R}}
\newcommand{\Z}{{\mathbb Z}}

\newcommand{\Q}{{\mathbb Q}}
\newcommand{\N}{{\mathbb N}}

\newcommand{\A}{{\mathcal A}}
\newcommand{\B}{{\mathcal B}}

\title{When  is the fractal uncertainty principle for discrete Cantor sets most uncertain?}

\author{Chun-Kit Lai}
\address{Department of Mathematics, San Francisco State University,
1600 Holloway Avenue, San Francisco, CA 94132.}
 \email{cklai@sfsu.edu}

\author{Ruxi Shi}

\address{Shanghai Center for Mathematical Sciences, Fudan University, 200438 Shanghai, China.}
 \email{ruxishi@fudan.edu.cn}

\begin{document}
\begin{abstract}
We give a necessary and sufficient condition to achieve the most uncertain exponent in the fractal uncertainty principle of discrete Cantor sets. The condition will be described as distributed spectral pairs, which is a generalization of the spectral pair studied in the spectral sets literature. We investigate distributed spectral pairs in some cyclic groups and some complete classifications are given. Finally, we also discuss the most uncertain case in the continuous setting. 
\end{abstract}
\date{\today}

\maketitle

\section{Introduction}
A {\bf fractal uncertainty principle} (FUP), roughly speaking, states that {\it no function can be localized both close to a fractal set in both positions and frequencies}. It has found striking applications in the control of eigenvalues and the spectral gap problem of operators that appear in quantum chaos and hyperbolic dynamics. Common FUP takes place in both continuous and discrete settings.  Readers are invited to consult Dyatlov's article \cite{Dyatlov2019} and the reference therein for its applications and some of the basic results in this field.



 The present article will mainly focus on the discrete setting and we are interested in the largest uncertainty exponents. Let us first review the  discrete version of the fractal uncertainty principle introduced by Dyatlov and Jin \cite{DyaJin}. 
Let $M\ge 3$ be an integer and $\Z_M$ be the cyclic group of $M$ elements. When we write $A\subset\Z_M$, it is possible that elements of $A$ may be larger than $M$.   We will denote by ${\mathcal F}_M$ the standard discrete Fourier transform matrix, i.e.
$$
{\mathcal F}_M = \frac{1}{\sqrt{M}} \left(e^{-2\pi i \frac{jk}{M}}\right)_{0\le j,k\le M-1}.
$$
Fix an alphabet set ${\mathcal A}\subset \{0,1,...,M-1\}$, we will always assume that $1<|\A|<M$ to avoid triviality ($|\A|$ denotes the cardinality of $\A$). In particular, 
$$
\delta := \frac{\log |\A|}{\log M}\in (0,1).
$$
Define the $k^{th}$ order iteration  by 
\begin{equation}\label{A_k}
{\mathcal A}_k  := \left\{ \sum_{j=0}^{k-1} a_j M^j: a_j\in\mathcal{A}\right\}  = {\mathcal A}+M{\mathcal A}+\cdots+M^{k-1}{\mathcal A}.
\end{equation}
Note that $\A_k\subset \{0,1,\cdots, M^k-1\}$ if $\A\subset \{0,1,\cdots M-1\}$.  $(M, \A)$ generates a self-similar Cantor set 
\begin{equation}\label{eq_CMA}
C = C(M, \A) = \bigcap_{k=1}^{\infty}\bigcup_{a\in\A_k} \left[\frac{a}{M^k}, \frac{a+1}{M^k}\right].
\end{equation}
It is well-known that $C(M,\A)$ is $\delta$-Ahlfors-David regular with Hausdorff dimension $\delta$. Consider two alphabet sets ${\mathcal A}$ and ${\mathcal B}$ where $|{\mathcal A}| = |{\mathcal B}|$.  We  denote ${\bf 1}_E$ the indicator function on the set $E\subset {\mathbb Z}_M$. It is now well-known (see e.g. \cite{Dyatlov2019}) that the operator ${{\mathbf 1}}_{{\mathcal A}_k} {\mathcal F}_{M^k}{\bf 1}_{{\mathcal B}_k}$ obeys the following bound:  
$$
\|{{\mathbf 1}}_{{\mathcal A}_k} {\mathcal F}_{M^k}{\bf 1}_{{\mathcal B}_k}\|_{{\ell^2\to\ell^2}}\lesssim M^{-k\max\{0,\frac12-\delta\}}.
$$
 Indeed, ${\mathbf 1}_{{\mathcal A}_k} {\mathcal F}_{M^k}{\bf 1}_{{\mathcal B}_k}$ is the submatrix of ${\mathcal F}_{M^k}$ by choosing columns corresponding to  ${\mathcal B}_k$ and rows corresponding to  ${\mathcal A}_k$ and $\|\cdot\|_{\ell^2\to\ell^2}$ denotes the operator norm of the matrix acting on the finite-dimensional vector space ${\mathbb C}^{|\A|^k}$ with Euclidean norm.   Dyatlov and Jin \cite{DyaJin} proved the following fractal uncertainty principle in the discrete Cantor set setting.

\begin{theorem}\label{th_DJ}
Let ${\mathcal A},{\mathcal B}\subset \{0,1,...,M-1\} $. Then there exists $\beta  =\beta (M,{\mathcal A},{\mathcal B})> \max \left(0,\frac12-\delta\right)$ such that 
\begin{equation}\label{FUP1}
\|{{\mathbf 1}}_{{\mathcal A}_k} {\mathcal F}_{M^k}{\bf 1}_{{\mathcal B}_k}\|_{\ell^2\to\ell^2} \lesssim M^{-k\beta}.
    \end{equation}
 \end{theorem}

\medskip

The largest exponent $\beta (M, {\mathcal A},{\mathcal B})$ for which (\ref{FUP1}) holds  will be  called  the {\bf uncertain exponent} for the discrete Cantor sets generated by $M,{\mathcal A}$ and ${\mathcal B}$.  The largest exponent could vary when we use different $\A,\B$ and it can be very close to no improvement (see e.g. \cite[Proposition 3.17]{DyaJin}). Eswarathasan and Han \cite{EH2023} investigated the best exponents for random alphabets, and with high probability the exponent is $\frac12-\frac34\delta$ for $\delta<\frac23$.

\subsection{The most uncertain exponent} In \cite{DyaJin}, it was observed that by taking the vector that supports at only one point of $Y$ or X (by considering $({\bf 1}_X{\mathcal F}_N{\bf 1}_Y)^{\ast}$), we have
$$
\|{\bf 1}_X{\mathcal F}_N{\bf 1}_Y\|_{\ell^2\to\ell^2}\ge \sqrt{\frac{\max\{|X|,|Y|\}}{N}}.
$$
Consequently, by taking $X = \A_k$, $Y  = \B_k$ and $N = M^k$,
$$
\beta (M, \A, \B)\le \frac{1-\delta}{2}. 
$$
That means that the largest (or the most uncertain) exponent can at most be $\frac{1-\delta}{2}$.   The main purpose of this paper is to study for what alphabets ${\mathcal A}$ and ${\mathcal B}$ with respect to the scale $M$,  $\beta (M, \A, \B)$ achieves the most uncertain exponent. This turns out to be related to the study of  spectral sets initiated first by Fuglede \cite{Fuglede1974}.

\medskip
\begin{definition}
     A {\bf spectral pair} in the cyclic group $\Z_M$ is a pair of subsets $(\A,\B)$  with $|\A| = |\B|$ such that the matrix ${\bf 1}_A{\mathcal F}_{M}{\bf 1}_B$ has mutually orthogonal columns. Equivalently, the matrix 
     $$
     \frac{1}{\sqrt{|\A|}}\left(e^{2\pi i \frac{ab}{M}}\right)_{a\in \A,b\in\B} = ~\sqrt{\frac{M}{|\A|}}\cdot{\bf 1}_A{\mathcal F}_{M}{\bf 1}_B
     $$
     is a unitary matrix on the complex vector space ${\mathbb C}^{|\A|}$. 
\end{definition}
 A spectral pair in $\Z_M$ will naturally attain the most uncertain exponent.  Indeed, this follows from an easy mathematical induction that $({\mathcal A}_k,{\mathcal B}_k)$ forms a spectral pair in $\Z_{M^k}$ for all $k\ge 1$. If we let
${\mathcal H}_k = {{\mathbf 1}}_{{\mathcal A}_k} {\mathcal F}_{M^k}{\bf 1}_{{\mathcal B}_k}$, then 
$${\mathcal H}_k^{\ast}{\mathcal  H}_k = \left(\frac{|{\mathcal A}|}{M}\right)^kI.$$
Hence, we obtain that
$$
\|{\bf 1}_A{\mathcal F}_{M}{\bf 1}_B\|_{\ell^2\to\ell^2} = \left(\frac{|{\mathcal A}|}{M}\right)^{k/2}
 = M^{-k\frac{1-\delta}{2}}$$
 This shows that a spectral pair induces immediately the most uncertain case.

The most uncertain exponent is connected with {\it Spectral set conjecture} or so-called {\it Fuglede's conjecture}: a Borel set is a spectral set in $\R^n$ if and only if it tiles $\R^n$ by translation. Recall that a Borel set $\Omega$ is a {\bf spectral set} in $\R^n$  if there exists a discrete set $\Gamma$ such that $\{e^{2\pi i \gamma x}\}_{\gamma\in \Gamma}$ forms a Hilbert basis of $L^2(\Omega)$.  Although the conjecture was disproved
eventually for $\R^d$ with $d\ge 3$, it remains still widely
open for general locally compact abelian groups, especially for abelian
groups in lower “dimension”. We only know that Fuglede's spectral
set conjecture holds on some specific groups, for example,  
$\Z_p\times \Z_p$ \cite{IMP2017}, $\Z_{p^n}$ and $\Z_{p^nq}$ with $n\ge 1$ \cite{MK17}, $p$-adic field $\Q_p$ \cite{FFS2017,FFLS19}, $\Z_{pqr}$ with $p,q,r$
diﬀerent primes \cite{Shi2019}, $\Z_{p^nq^m}$ with $m$ small \cite{M2022}, $\Z_{p^2} \times \Z_p$ \cite{S2020}, $\Z_{pqrs}$ with $p,q,r,s$
diﬀerent primes \cite{KMSV2022} and $\Z_{p^2} \times \Z_{p^n}$ \cite{Z2023}.

 \subsection{Main results} One of the main questions is to determine what condition for which $(\A,\B)$ achieves the most uncertain exponent. One would naturally ask if spectral pairs are the only case that the most uncertain exponent happens.  For a given finite subset $E\subset \R$, let $\nu_E$ denote the normalized probability measure supported on $E$. Its Fourier transform is defined by 
$$
\widehat{{\nu}_E} (\xi) = \frac{1}{|E|}\sum_{a\in E} e^{-2\pi i \xi a}.
$$
Notice that if $(\A,\B)$ is a spectral pair in $\Z_M$, then  
$$
\widehat{{\nu}_{\mathcal A}} \left(\frac{b-b'}{M}\right) = 0,
$$
for all $b\ne b'\in \B$. Using our notation\footnote{We normalize so that $\nu_E$ is a probability measure. The normalization constant in \cite{DyaJin} was $M^{-1/2}$, which is different from ours.},  a partial converse concerning the most uncertain exponent was obtained by Dyatlov and Jin \cite[Proposition 3.16]{DyaJin}.

\begin{proposition}\label{Prop_DJ}
Suppose that $\beta = \frac{1-\delta}{2}$. Then for all $b\ne b'$, $b,b'\in{\mathcal B}$, we have
\begin{equation}\label{eqDJ}
  \widehat{{\nu}_{\mathcal A}} \left(\frac{b-b'}{M}\right)\widehat{{\nu}_{\mathcal A}} \left(\frac{b-b'}{M^2}\right) = 0.  
\end{equation}
\end{proposition}

In this paper, we provide a complete characterization for which $(M,\A,\B)$ has its most uncertain exponent. 

\begin{theorem}\label{maintheorem1}
Suppose that ${\mathcal A},{\mathcal B}\subset \{0,1,\cdots M-1\}$ and we form the discrete Cantor sets ${\mathcal A}_k,{\mathcal B}_k$ as before. Then the following are equivalent. 
\begin{enumerate}
\item $\beta (M,{\mathcal A},{\mathcal B}) = \frac{1-\delta}{2}$
\item for all $b_1\ne b_1'\in{\mathcal B}$, 
$${\widehat{{\nu}_{\mathcal A}} \left(\frac{b_1-b_1'}{M}\right)} = 0 ~~~\mbox{or}~~\\\ ~ {\widehat{{\nu}_{\mathcal A}} \left(\frac{b_1-b_1'}{M^2}+ \frac{b_2-b_2'}{M}\right)} = 0  \ \textrm{for all} \  b_2,b_2'\in \B
$$.
\end{enumerate}
\end{theorem}

We now say that $(\A,\B)$ is a {\bf distributed spectral pair} in $\Z_M$ (not necessarily in $ \{0,1,\cdots M-1\}$) if it satisfies the condition (2) in Theorem \ref{maintheorem1}.  From condition (2) above, if  $b_1-b_1'$ is not zero in the first factor, its zero must be distributed to all translates by $(b_2-b_2')/M$ in the second factor. This theorem offers a readily checkable criterion for which $\beta = \frac{1-\delta}{2}$ if $\A,\B$ are given. A spectral pair must be a distributed spectral pair. However, the converse is not true. Indeed, a simple example is given below. 

\begin{example}\label{example1.5}
Let ${\mathcal A}=\{0, 8\}$ and ${\mathcal B}=\{0,9\}$. Then $({\mathcal A}, {\mathcal B})$ is a spectral pair in $\Z_{12^2}$. However, $(\A,\B)$ is not a spectral pair in $\Z_{12}$. Indeed, as $a (b-b')\in 12 \Z$ for all $a\ne0$ and $b,b'\in\B$, it satisfies that $\widehat{\nu_{\A}}(\frac{b-b'}{M^2}+ \frac{b_2-b_2'}{M}) = 0$ for all $b\ne b'$ and $b_2,b_2'\in\B$. 
\end{example}

The following example shows that distributed spectral pairs are strictly stronger than the condition \eqref{eqDJ}. Therefore, the condition provided by Dyatlov and Jin in Proposition \ref{Prop_DJ} is not sufficient. 

\begin{example}
Let $\A = \{0,1,9,10\}$ and $\B = \{0,2,8,10\}$ in $\Z_{12}$. Then
$$
B - B = \{0,\pm2,\pm 6,\pm 8,\pm10\}.
$$
From a direct check, 
$$
\widehat{\nu_{\A}} \left(\frac{\pm2}{12}\right) = \widehat{\nu_{\A}} \left(\frac{\pm6}{12}\right) = \widehat{\nu_{\A}} \left(\frac{\pm8}{(12)^2}\right) = \widehat{\nu_{\A}} \left(\frac{\pm10}{12}\right) = 0.
$$
Hence, condition (\ref{eqDJ}) holds. However, $\widehat{\nu_{\A}} \left(\frac{8}{12}\right)\ne 0$ and 
$$
\widehat{\nu_{\A}} \left(\frac{8}{12^2}+\frac{2}{12}\right) = \widehat{\nu_{\A}} \left(\frac{2}{3^2}\right) \ne 0.
$$
Therefore, $(\A,\B)$ is not a distributed spectral pair. By Theorem \ref{maintheorem1}, $\beta (M,\A,\B)<\frac{1-\delta}{2}$. 
\end{example}

\subsection{Overview of the paper.} The proof of Theorem \ref{maintheorem1} will be given in Section \ref{proof of thm}. The sufficiency of the theorem will be proved via an inductive construction of the Hermitian matrices. We notice that the sufficiency part does not require $\A,\B$ to be a subset of $\{0,1,\cdots M-1\}$.  The necessity part of the theorem is  motivated by the proof of Proposition \ref{Prop_DJ}, but the estimates are reorganized using the perspective of the associated self-similar measures. In this proof, the fact that $\A,\B$ is a subset of $\{0,1,\cdots M-1\}$ plays an important role. 

It is now an interesting question to ask if we can have a further classification of distributed spectral pairs. However, it does not seem to be an easy question. Section \ref{Spectral pair} and Section \ref{classification} will be devoted to such study. 

In Section  \ref{Spectral pair}, we will show that  spectral pairs in $\Z_{M^2}$ with $\A,\B\in\{0,1,\cdots M-1\}$ are indeed distributed spectral pairs and they contain only two elements. However, these pairs exist in many cyclic groups depending on the relative size of the prime factors.   In Section \ref{classification}, we give a complete  classification of the distributed spectral pair in $\Z_M$ for  $M = p^k$, $M = pq$, and $M=p^2q$ respectively where $p,q$ are distinct prime numbers. We will show that, among these groups, the only non-trivial distributed spectral pair is the one from Example \ref{example1.5}. 

 FUP for discrete Cantor sets are also studied in higher dimension \cite{C2022}. It is also straightforward to see that a distributed spectral pair can attain the largest exponents if we consider the self-similar setting (i.e. $M$ is replaced by a diagonal matrix of all entries equal to $M$). We believe that with a similar strategy of our proof, we can also obtain similar results in higher dimensions. 

 In Section \ref{discussion}, we will wrap up some problems and observations in the discrete settings. Finally, we will also discuss the most uncertain exponent for the continuous FUP. We anticipate that the limiting sets of the distributed spectral pairs will attain the most uncertain exponents as well, though we do not yet have a concrete proof right now.

\section{Proof of Theorem \ref{maintheorem1}} \label{proof of thm}

With the notion of distributed spectral pair, we rewrite Theorem \ref{maintheorem1} as follows.

\begin{theorem} (=Theorem \ref{maintheorem1})\label{maintheorem11}
\begin{enumerate}
\item  Let ${\mathcal A},{\mathcal B}\subset \Z_M$ and suppose that $({\mathcal A}, {\mathcal B})$ is a distributed spectral pair.  Then it achieves the most uncertain exponent, i.e.  $\beta (M,{\mathcal A},{\mathcal B}) = \frac{1-\delta}{2}$. 
\item Conversely, if ${\mathcal A},{\mathcal B}\subset \{0,1,\cdots M-1\}$ and $\beta (M,{\mathcal A},{\mathcal B}) = \frac{1-\delta}{2}$. Then $({\mathcal A}, {\mathcal B})$ is a distributed spectral pair. 
\end{enumerate}
\end{theorem}

In this section, we prove Theorem \ref{maintheorem11}.
We first notice a translational invariance property for our operator. If we let $\widetilde{B} = B-b_0$ and $\widetilde{A} = A-a_0$ and define the corresponding $\widetilde{B}_k$ and $\widetilde{A}_k$, then 
$$
\|{\bf 1}_{\widetilde{A}_k}{\mathcal F}_{M^k}{\bf 1}_{\widetilde{B}_k}\|_{\ell^2\to\ell^2} = \|{\bf 1}_{A_k}{\mathcal F}_{M^k}{\bf 1}_{{B}_k}\|_{\ell^2\to\ell^2}.
$$
Hence, we will assume that $0\in \A\cap \B$. We notice a simple linear algebra fact that if $\lambda\ne 0$, then $\lambda$ is an eigenvalue for the product $A^{\ast}A$ if and only if $\lambda$ is an eigenvalue for $AA^{\ast}$. Hence, this implies that 
\begin{equation}\label{eq-dual}
\beta(M,\A,\B) = \beta (M,\B,\A).
\end{equation}

To prove the theorem, we will need to set up some notations. The matrix of our interest will be denoted by 
$$
{\mathbf F}_k = ({{\mathbf 1}}_{{\mathcal A}_k} {\mathcal F}_{M^k}{\bf 1}_{{\mathcal B}_k})^{\ast}({{\mathbf 1}}_{{\mathcal A}_k} {\mathcal F}_{M^k}{\bf 1}_{{\mathcal B}_k}).
$$
Note that $\|({{\mathbf 1}}_{{\mathcal A}_k} {\mathcal F}_{M^k}{\bf 1}_{{\mathcal B}_k})\|_{\ell^2\to\ell^2}$ is equal to  the  square root of largest eigenvalue for the matrix ${\mathbf F}_k$, which is equal to the following $|{\mathcal B}|^k\times|{\mathcal B}|^k$ matrix: 
$${\mathbf F}_k =  \left[ \frac1{M^k}\sum_{{\bf a}\in{\mathcal A}_k} e^{2\pi i \frac{{\bf a} ({\bf b}-{\bf b}')}{M^k}}\right]_{{\bf b},{\bf b}'\in {\mathcal B}_k} =  \frac{|\A|^k}{M^k}\left[\widehat{{\nu}_{\A_k}} \left( \frac{{\bf b}-{\bf b}'}{M^k}\right) \right]_{{\bf b},{\bf b}'\in {\mathcal B}_k}.
$$
Moreover, using the equation (\ref{A_k}) for ${\mathcal A}_k$, 
\begin{equation}\label{FT_A_K}
\widehat{{\nu }_{{\mathcal A}_k}} \left(\frac{\xi}{M^k}\right) = \prod_{j=1}^{k} \widehat{{\nu}_{{\mathcal A}}} \left(\frac{\xi}{M^j}\right) = \widehat{\nu_{\A}} \left(\frac{\xi}{M}\right)\cdot \widehat{\nu_{\A_{k-1}}} \left(\frac{\xi}{M^k}\right)
\end{equation}
We now decompose ${\bf b},{\bf b}'\in\B_k$ for $k\ge 3$ into 
\begin{equation}\label{equation_b_decompose}
{\bf b} = b_1+  M b_2+ M^2 {\bf b}_{k-2}, \ {\bf b}' = b_1'+  M b_2'+ M^2 {\bf b}_{k-2}'
\end{equation}
where $b_1,b_1',b_2,b_2'\in \B$ and ${\bf b}_{k-2},{\bf b}_{k-2}'\in\B_{k-2}$. Define also ${\bf b}_{k-1} = b_2+M{\bf b}_{k-2}$ and ${\bf b}'_{k-1} = b'_2+M{\bf b}'_{k-2}$. Plugging ${\bf b}-{\bf b}'$ into (\ref{FT_A_K}) and note that $\widehat{\nu_A}$ is integer periodic, we have
\begin{equation}\label{equation_iteration}
\begin{aligned}
\widehat{{\nu}_{\A_k}}\left(\frac{{\bf b}-{\bf b}'}{M^k}\right) = &  \widehat{\nu}_{\A} \left(\frac{b_1-b_1'}{M}\right)\widehat{\nu}_{\A} \left(\frac{b_1-b_1'}{M^2}+\frac{b_2-b_2'}{M}\right) \cdot \widehat{\nu}_{\A_{k-2}} \left(\frac{{\bf b}-{\bf b}'}{M^k}\right).
\end{aligned}
\end{equation}
\medskip

\subsection{Proof of Theorem \ref{maintheorem11}(1)} We now prove that  a distributed spectral pair must be the most uncertain.  We will show by induction that for all $k\ge 3$, we can decompose ${\mathbf F}_k$ into $|\B|\times |\B|$ block matrix with each block has size $|\B|^{k-1}\times |\B|^{k-1}$ as follows:
\begin{equation}\label{equation_block1}
{\bf F}_k = \frac{|\A|}{M}\left[\begin{array}{cccc} {\mathbf F}_{k-1} & O &....& O \\ O& {\mathbf F}_{k-1}  &....& O\\  & & \ddots& \\ O & O & \cdots & {\mathbf F}_{k-1} \end{array}\right]
\end{equation}
($O$ denotes the zero matrix). To see this, according to  the decomposition (\ref{equation_b_decompose}), we decompose the matrix into blocks according to $b_1$ and $b_1'$. If  $b_1\ne b_1'$,  applying (\ref{equation_iteration}) and our assumption implies that the first two factors in (\ref{equation_iteration}) is always zero. Hence, we must have an $|\B|^{k-1}\times |\B|^{k-1}$ zero matrix. If $b_1=b_1'$, then the first term in (\ref{equation_iteration}) becomes $1$, so we have
$$
\begin{aligned}
\widehat{{\nu}_{\A_k}}\left(\frac{{\bf b}-{\bf b}'}{M^k}\right) = &\widehat{{\nu}_{\A}}\left(\frac{b_2-b_2'}{M}\right) \widehat{\nu}_{\A_{k-2}} \left(\frac{b_2-b_2'+M({\bf b}_{k-2}-{\bf b}_{k-2}')}{M^{k-1}}\right)\\
=& \widehat{{\nu}_{\A_{k-1}}}\left(\frac{{\bf b}_{k-1}-{\bf b}_{k-1}'}{M^{k-1}}\right).
\end{aligned}
$$ 
which means that the matrix reduces to the $\frac{|\A|}{M}{\mathbf F}_{k-1}$. This justifies (\ref{equation_block1}).

\medskip

Because of (\ref{equation_block1}) and the fact that the operator norm $\|A\|_{\ell^2\to\ell^2}$ is equal to the square root
of the largest eigenvalue of $A^{\ast} A$, we have 
$$
\|{{\mathbf 1}}_{{\mathcal A}_k} {\mathcal F}_{M^k}{\bf 1}_{{\mathcal B}_k}\|_{\ell^2\to\ell^2} = \left(\frac{|A|}{M}\right)^{\frac{k-2}{2}} \|{{\mathbf 1}}_{{\mathcal A}_2} {\mathcal F}_{M^2}{\bf 1}_{{\mathcal B}_2}\|_{\ell^2\to\ell^2} =  C \cdot M^{-k\cdot \frac{1-\delta}{2}}
$$
for some constant $C>0$ independent of $k$. Hence, (1) holds. 

\subsection{Proof of Theorem \ref{maintheorem11}(2).} We now prove that  a distributed spectral pair must be necessary when $\A,\B\subset \{0,1,\cdots M-1\}$. We suppose that $(\A,\B)$ does not form a distributed spectral pair.  Then there exists $b_1\ne b_1'$  and $b_2, b_2'\in\B$ such that 
\begin{equation}\label{equation_contradiction}
\widehat{{\nu}}_{\A} \left(\frac{b_1-b_1'}{M}\right)\ne 0 \ \mbox{and} \ \widehat{{\nu}}_{\A} \left(\frac{b_1-b_1'}{M^2}+\frac{b_2-b_2'}{M}\right)\ne 0. 
\end{equation}

\medskip

We now show that $\beta(M,\B,\A)<\frac{1-\delta}{2}$, then our proof will be complete. The idea is adapted from the proof of Proposition \ref{Prop_DJ}. Let $k$ be an integer  and let $r_{2k} = \|{\bf 1}_{\B_{2k}}{\mathcal F}_{M^{2k}}{\bf 1}_{\A_{2k}}\|_{\ell^2\to\ell^2}$. 
Then writing ${\bf u}$ be the the constant vector with all entries 1, we have 
\begin{equation}\label{eqr_k}
r_{2k}^2 \ge \frac{\|{\bf 1}_{\B_{2k}}{\mathcal F}_{M^{2k}} {\bf 1}_{\A_{2k}}{\bf u}\|^2}{\|{\bf 
 u}\|^2}  = \frac{1}{\|{\bf u}\|^2}\sum_{{\bf b}\in \B_{2k}} \left|\sum_{a\in \A_{2k}} \frac{1}{\sqrt{M^{2k}}} e^{-2\pi i \frac{a{\bf b}}{M^{2k}}}\right|^2.
\end{equation}
Then $\|{\bf u}\|^2 = |\A|^{2k}$ and plugging back into (\ref{eqr_k}) and using (\ref{FT_A_K}), we have
\begin{equation}\label{eqr_k-1}
r_{2k}^2\ge \left(\frac{|\A|}{M}\right)^{2k}\sum_{{\bf b}\in \B_{2k}} \left|\widehat{\nu}_{{\A}_{2k}} \left(\frac{{\bf b}}{M^{2k}}\right)\right|^2 = \left(\frac{|\A|}{M}\right)^{2k}\sum_{{\bf b}\in \B_{2k}} \left|\prod_{j=1}^{2k}\widehat{\nu}_{{\A}} \left(\frac{{\bf b}}{M^{j}}\right)\right|^2.
\end{equation}
\subsubsection{\bf Associated self-similar measures.} In estimating the product above, it will be useful to define the associated  self-similar measure $\mu(M,\A)$ supported on the Cantor set $C(M,\A)$.  This measure is an infinite convolution of discrete measures 
\begin{equation}\label{eq-self-similar}
\mu =\mu(M,\A) = \nu_{\frac{\A}{M}}\ast\nu_{\frac{\A}{M^2}}\ast...    
\end{equation}
$\mu$ defines uniquely a self-similar probability measure supported on the Cantor set $C(M,\A)$ in (\ref{eq_CMA}) satisfying the following invariant identity:
$$
\mu(E) = \sum_{a\in \A} \frac{1}{|\A|} \mu (M E-a), \ \forall E \ \mbox{Borel}.
$$ 
Details of this measure can be found in \cite{Fal97}. We can compute its Fourier transform from the infinite convolution product:
$$
\widehat{\mu}(\xi) = \prod_{j=1}^{\infty} \widehat{\nu_{\A}}\left(M^{-j}\xi\right).
$$
It defines an analytic function that is uniformly continuous on $\R$.  In our situation, we will factorize consecutive terms of the self-similar measures so that 
$$
\widehat{\mu}(\xi) = \prod_{j=1}^{\infty} \widehat{\nu_{\A_2}}\left(M^{-2j}\xi\right)
$$
Let $\xi_0 = (b_1-b_1')+M(b_2-b_2')$ 
where $b_1,b_1'b_2,b_2'$ are defined as in (\ref{equation_contradiction}). In the following, we  define 
$$
\xi_0 = \left\{\begin{array}{ll}
b_1-b_1'& \ \mbox{if} ~ \widehat{\nu_A} \left(\frac{b_1-b_1'}{M^2}\right)\ne0, \\
(b_1-b_1')+M(b_2-b_2')& \mbox{if} ~ \widehat{\nu_A} \left(\frac{b_1-b_1'}{M^2}\right)=0. 
\end{array}\right.
$$

\begin{lemma}\label{lem-non-zero}
    $\widehat{\nu_{\A_2}}(M^{-2j}\xi_0)\ne 0$ if $j>1$.
\end{lemma}

\begin{proof}
    Note that $-M^2<\xi_0<M^2$ since $\B\subset\{0,1,\cdots, M-1\}$. If $j\ge 3$, then $e^{2\pi i aM^{-2j}\xi_0}$ lies on the upper half-plane for all $a\in \A_2$, so the $\widehat{\nu_{\A_2}}(M^{-2j}\xi_0)$ cannot be zero. In the same reason, when $j=2$ and $\widehat{\nu_A} \left(\frac{b_1-b_1'}{M^2}\right)\ne0$, then $e^{2\pi i aM^{-4}(b_1-b_1')}$ lies on the upper half-plane for all $a\in\A_2$. $\widehat{\nu_{\A_2}}(M^{-4}\xi_0)\ne0$ 
    
     It remains to show the lemma for the case  $\widehat{\nu_\A} \left(\frac{b_1-b_1'}{M^2}\right)=0$. we need to claim that  $\widehat{\nu_{\A_2}}(M^{-4}\xi_0) = \widehat{\nu_{\A_2}}(M^{-4}\xi_0)\ne0$. We will argue by contradiction. Suppose that    $\widehat{\nu_{\A_2}}(M^{-4}\xi_0)=0$, then using $\A_2 = \A+M\A$. we have
     $$
     \widehat{\nu_{\A}} (M^{-4}\xi_0)\cdot\widehat{\nu_{\A}} (M^{-3}\xi_0) =0.
     $$
     We must have $\widehat{\nu_{\A}} (M^{-3}\xi_0) =0$ since $ \widehat{\nu_{\A}} (M^{-4}\xi_0)\ne0$ using the fact that $e^{2\pi i aM^{-4}\xi_0}$ are on the upper half-plane. Let 
     $$
     P_{\A}(z) = \sum_{a\in A}z^a, \ \mbox{and} \ f(x) = P_{\A}(e^{2\pi i x}).$$
     Then $f(x)$ has at least two distinct real roots in $(0,1/M)$, namely $x_1 = (b_1-b_1')/M^2$ and $M^{-3}\xi_0$. Therefore, counting conjugate roots,  $P_{\A}$ has at least four roots in the sector $|\arg(z)|< 2\pi/M$. As degree of $P_\A$ is at most $M-1$, this is impossible if $M\le 4$.

     When $M\ge 5$, we apply a result of Obrechkoff  \cite{Obrechkoff1923Laguerre} (see also \cite{Obrechkoff}), which says that  {\it a degree $d$ polynomial with nonnegative coefficients has at most $2d\alpha/\pi$ many roots in the sector where $|\arg(z)|< \alpha$}. In the sector of an opening angle $\alpha = 2\pi/M$, $P_{\A}(z)$ has at most $4(M-1)/M<4$ many roots. This is a contradiction since we know $P_{\A}$ has at least four roots.\footnote{The last part of the proof was suggested by ChatGPT}. 
\end{proof}

\medskip

 Lemma \ref{lem-non-zero} together with our assumption (\ref{equation_contradiction}) implies that $\widehat{\mu}(\xi_0)\ne 0$. By continuity, there exists $\epsilon_0 = \frac12|\widehat{\mu}(\xi_0)|>0$ and $L>0$ such that 
$$
|\widehat{\mu}(\xi)| \ge \epsilon_0>0, \ \forall \xi\in [\xi_0(1-M^{-2L}), \xi_0(1+M^{-2L})].
$$
We will choose a large $L$ in the later argument.  Fixing $k> 2L$ and  $q>0$ such that $q\le \frac{k}{2L}$, consider the following subset of $\B_{2k}$:
$$
\Omega_{q,k,L} = \left\{\sum_{p=1}^{q} \xi_0 M^{2s_p}: 0\le s_1,...,s_q<k, \ s_{p+1}-s_p\ge L+1, \forall 1\le p<q \right\}.
$$
 Note that ${\bf b}\in \Omega_{q,k,L}$ if and only if we can write 
$$
{\bf b} = \sum_{p=1}^{q} \xi_0 M^{2s_p}.
$$
Hence, if we define  $s_0 = 0$, then 
$$
\begin{aligned}
\prod_{j=1}^{2k}\widehat{\nu}_{{\A}} \left(\frac{{\bf b}}{M^{j}}\right) = &\prod_{p=0}^q \prod_{j={s_p+1}}^{s_{p+1}} \widehat{\nu}_{{\A_2}} \left(\frac{{\bf b}}{M^{2j}}\right) \\
\ge & \prod_{p=0}^q \prod_{j={1}}^{s_{p+1}-s_p}\widehat{\nu}_{\A_2} \left(\frac{\xi_0}{M^{2j}}\left(1+\frac{1}{M^{2(s_p-s_{p-1})}}+...+\frac{1}{M^{2(s_p-s_{1})}}
\right)\right)\\
\ge & \prod_{p=0}^q \widehat{\mu} \left(\xi_0\left(1+\frac{1}{M^{2(s_p-s_{p-1})}}+...+\frac{1}{M^{2(s_p-s_{1})}}\right) \right)\\
\ge & \epsilon_0^q.
\end{aligned}
$$
Thus, (\ref{eqr_k-1}) implies that 
\begin{equation}\label{eq_omega}
r_{2k}^2\ge \left(\frac{|\A|}{M}\right)^{2k} \sum_{q = 1}^{\lfloor\frac{k}{2L}\rfloor}\epsilon_0^q\cdot |\Omega_{q,k,L}|. 
\end{equation}

\subsubsection{\bf Cardinality of $\Omega_{q,k,L}$.} The cardinality of $\Omega_{q,k,L}$ was computed without proof in Dyatlov and Jin's paper that $\frac{(k- (q-1)L)!}{q!(k-qL)!}$. However, we believe the formula is not correct. When $q=2$, we can count explicitly the number of elements in $\Omega_2 = \Omega_{2,k,L}$ by listing them out 
$$
\Omega_2 = \{\xi_0M^{2a}+\xi_0M^{2b}: 0\le a< k-L-1, \ a+L+1\le b\le k-1\}
$$
Then if $a$ is fixed, there are $k-1-a-L$ choices of $b$, so $|\Omega_2| = (k-L-1)+ (k-L-2)+...+1= \frac{(k-L-1)(k-L)}{2}$, which is not their formula.  We have the following proposition computing the exact formula.

\begin{proposition}\label{Prop_Omega}
    $$
    |\Omega_{q,k,L}| = \binom{k-(q-1)L}{q}.
    $$
\end{proposition}

\begin{proof}We would like to thank the referee for suggesting a simpler proof. Under the mapping $t_j = s_j - (j-1)L$, the cardinality of $|\Omega_{q,k,L}|$ is equal to the cardinality of the set
$$
\{(t_1,\cdots , t_q) : 0 \le  t_1 <t_2 <\cdots < t_q \le  k -1 - (q - 1)L\},
$$
which is equal to the number of $q$-elements subsets of $\{0,1\cdots, k-1-(q-1)L\}$. This is our desired formula. 
\end{proof}

\subsubsection{\bf Completion of the proof.}  Coming back to (\ref{eq_omega}), we note that $q\le k/L\le k/4$ and take $L\ge 4$, so Proposition \ref{Prop_Omega} implies that we can bound
$$
\begin{aligned}
|\Omega_{q,k,L}| = &\frac{(k-(q-1)L)(k-(q-1)L-1)...(k-(q-1)L-q+1))}{q!}\\
\ge &\frac{1}{q!}\left(\frac{q+1}{2q}k-\frac{k}{4}\right)^q \ge \frac{1}{q!}\left(\frac{k}{4}\right)^q,\\   
\end{aligned}
$$
for $q\ge 2$. Hence, with Stirling approximation, up to some constant $c>0$,
$$
\begin{aligned}
r_{2k}^2\ge &  \left(\frac{|\A|}{M}\right)^{2k} \sum_{q=2}^{\lfloor\frac{k}{2L}\rfloor} \left(\frac{\epsilon_0}{4}\right)^q \frac{k^q}{q!}\\
\ge &  c M^{-(1-\delta)2k} \left(\frac{\epsilon_0}{4}\right)^{k/2L} \frac{k^{k/2L}}{\left(\frac{k}{2L}\right)^{1/2} \left(\frac{k}{2L}\right)^{k/2L}e^{-k/2L}} \ \mbox{(keeping only the last term)}  \\
\ge & c M^{-(1-\delta)2k} \frac{\varrho^{2k}}{k^{1/2}}, 
\end{aligned}
$$
where  $\varrho = (2\epsilon_0Le/4)^{1/4L}>1$ if we choose $L$ sufficiently large. Consequently, $\beta \le \frac{1-\delta}{2}-\log_M\rho<\frac{1-\delta}{2}$, which completes the proof. 
\medskip

\section{Spectral pair in $\Z_{M^2}$}\label{Spectral pair}

In this section, we investigate spectral pairs in $\Z_{M^2}$ when they are contained in $\{0,1,...,M-1\}$.

\begin{proposition}\label{prop_spectral_M2}
Suppose that $(\A,\B)$ is a spectral pair in $\Z_{M^2}$ and $\A,\B\subset \{0,1,\cdots$ $,M-1\}$. Then $|\A |= |\B|= 2$. 
\end{proposition}

\begin{proof}
Note that if $b<\frac{M}{2}$, then $ab/M^2<1/2$. Thus, $\sum_{a\in\A}e^{2\pi i ab/M^2}\ne 0$ since all complex exponentials are on the upper half-plane. If $|\B|>2$ and it forms a spectral pair, then all distinct $b,b'\in \B$ must satisfy $|b-b'|\ge M/2$. But this is not possible since $\B\subset \{0,1,...,M-1\}$. This shows that $|\A |= |\B|= 2$. 
\end{proof}
\medskip
\begin{theorem}\label{theorem-spectral-m2}
Let $M\ge 2$ and let $\A,\B\subset \{0,1,...,M-1\}$ with $0\in \A\cap \B$. Then the following are equivalent. 
\begin{enumerate}
    \item $(\A,\B)$ is a spectral pair in $\Z_{M^2}$.
    \item $\A = \{0,a\}$ and $\B = \{0,b\}$, where 
    $$
    2ab = M^2, \ \mbox{and}  \  0< a, b<M.
$$
\end{enumerate}
Moreover,  $(\A,\B)$ is a distributed spectral pair in $\Z_M$. 
\end{theorem}
\begin{proof}
Suppose that $(\A,\B)$ is a spectral pair in $\Z_{M^2}$ and $\A,\B\subset \{0,1,...,M-1\}$. By Proposition \ref{prop_spectral_M2}, we know that  $|\A |= |\B|= 2$. Hence, we can write $\A = \{0,a\}$ and $\B = \{0,b\}$ for some $a,b\in \{0,1,\cdots M-1\}$.  Note that $(\A,\B)$ is a spectral pair in $\Z_{M^2}$ if and only if the matrix
$$
\begin{bmatrix}
    1 & 1 \\ 1 & e^{2\pi i ab/M^2}\\
\end{bmatrix}
$$
is an orthogonal matrix. Hence, the last exponential must be $-1$, implying that  $ab/M^2 = k/2$ where $k$ is an odd integer. Since $ab<M^2$, $k$ must be equal to 1. This shows that  $2ab  = M^2$. The rest of the conclusion is from our assumption.  Therefore, (2) holds. Conversely,  if (2) holds, the matrix above is orthogonal as the last entry is $-1$. The first part of the proof is complete. 

\medskip

It remains to show that $(\A,\B)$ is a distributed spectral pair in $\Z_{M}$. As $2ab  = M^2$, $M$ must be even. We now write $M = 2M'$ for some integer $M'$. Hence, $ab = M'\cdot M$. In particular, this implies that $ab\in M{\mathbb Z}$. This forces that 
$$
\widehat{\nu_{\A}}\left(\frac{b}{M^2}+ \frac{b}{M}\right) = \widehat{\nu_{\A}}\left(\frac{b}{M^2}\right) = 0. 
$$
As there are only two elements in $\B$, the above is the only case we need to check.  $(\A,\B)$ is thus a distributed spectral pair in $\Z_{M}$. 
\end{proof}

Theorem \ref{theorem-spectral-m2} gives a complete classification of $(\A,\B)$ is spectral in $\Z_{M^2}$ with $\A,\B\subset \{0,1,\cdots, M-1\}$. The following proposition examines more deeply how these pairs can be constructed. Not all groups have such pairs and it depends on the relative size of the prime factors.

\begin{proposition}\label{prop_prime-Char}
    
Let $M\ge 2$ be an integer. Let $\A = \{0,a\}$ and $\B = \{0,b\}$ with $a,b<M$ be a spectral pair in $\Z_{M^2}$. We have

\begin{enumerate}
    \item Suppose that $M = 2^{\alpha}q$ where $q$ is an odd prime. Then $(\A,\B)$ forming a spectral pair in $\Z_{M^2}$ if and only if $a = 2^{u}q^2$ and $b = 2^{2\alpha-u-1}$ and $2^{\alpha-u-1}<q<2^{\alpha-u}$ for some $0\le u<\alpha$. 
  
    \item If $M = 2p_1\cdots p_k$, where $p_1<\cdots<p_k$ are distinct odd primes and $a$ is even. Let 
    $$
    S_0 = \{1\le i\le k: p_i^2| a\}, \     S_1 = \{1\le i\le k: p_i| a \ \mbox{and} \ p_i|b\},  \ S_2 = \{1\le i\le k: p_i^2|b \}.
    $$
    
    Then $(\A,\B)$ forming a spectral pair in $\Z_{M^2}$ if and only if  
    \begin{equation}\label{eqab}
    a = 2\cdot\left( \prod_{i\in S_0}p_i^2\right)\cdot \left(\prod_{i\in S_1}p_i\right), \ b =\left( \prod_{i\in S_1}p_i\right)\cdot \left(\prod_{i\in S_2}p_i^2\right)\end{equation} 
    and 
\begin{equation}\label{eq_ab product}
     \prod_{i\in S_0} p_i  <\prod_{i\in S_2} p_i <  2\cdot \prod_{i\in S_0} p_i.
\end{equation}
\item Suppose that $M = 2pq$ where $p<q$ are distinct odd prime numbers and $a$ is even.  Then $(\A,\B)$ forming a spectral pair in $\Z_{M^2}$ if and only if $a=2p^2$ and $b = q^2$ and $q<2p$. 
\end{enumerate}

\end{proposition}

\begin{proof}

(1) From Theorem \ref{theorem-spectral-m2}, we know that $0< a,b< 2^{\alpha}q$ and $ab = 2^{2\alpha-1}q^2.$ We first observe that $q^2$ must divide $a$ or $b$. If not, then we can write $a = 2^uq$ and $b=2^v q$ where $u+v = 2\alpha-1$. On the other hand, $a<M$ implies that $u<\alpha$, so is $v$. Hence, $u+v\le2(\alpha-1)<2\alpha-1$. Therefore, such decomposition is not possible. 

\medskip

Without loss of generality, we write $a = 2^uq^2$ and $b = 2^v$ where $u+v = 2\alpha-1$. The condition that $a<M$ and $b<M$ holds if and only if $2^{v-\alpha}<q<2^{\alpha-u}$. But $\alpha-u = v-\alpha+1$. The condition is equivalent to 
$2^{\alpha-u-1}<q<2^{\alpha-u}$, which is exactly the desired statement.

\medskip

\noindent (2) From Theorem \ref{theorem-spectral-m2}, we know that $ab = M^2/2 = 2p_1^2\cdots p_k^2$, so $S_0\cup S_1\cup S_2 = \{1,\cdots k\}$ and the union is disjoint.  With $a$ being even, the definition of $S_i$ implies that $a,b$ must have the form (\ref{eqab}). Finally, $a<M$ and $b<M$ implies (\ref{eq_ab product}).

    \medskip

 \noindent (3) If $M = 2pq$,  we notice that $S_1 = \emptyset$. Otherwise, $a = 2p^2q$ or $2pq^2$ which are both larger than $M$. Hence, $a = 2p^2$ or $a= 2q^2$ or $b = p^2q$ or $pq^2$. Clearly, $2q^2>2pq$, so it is only possible that $a = 2p^2$ and thus $b = q^2$. Finally, (\ref{eq_ab product}) holds if and only if $p<q<2p$. 

\end{proof}

From the above proposition, we see that the following are spectral pairs in $\Z_{M^2}$ and they are also distributed spectral pairs. 

\begin{enumerate}
    \item $\A = \{0,2\cdot 3^2\}$ and $\B = \{0, 5^2\}$ in $\Z_{(2\cdot 3\cdot 5)^2}$. 
    \item $\A = \{0, 5^2\}$ and $\B = \{0, 2^5\}$ in $\Z_{(2^3\cdot 5)^2}$.
    \item $\A = \{0, 2\cdot 11\cdot 17\}$, $\B = \{0, 2\cdot 13\cdot 19\}$ in $\Z_{M^2}$, where $M = 2\cdot 11\cdot 13\cdot 17\cdot 19$.
    \item On the other hand, there is no spectral pair $(\A,\B)$  in $\Z_{(2^2\cdot q)^2}$ with $\A,\B\subset \{0,1,\cdots,M-1\}$ if $q>3$. This follows from Proposition \ref{prop_prime-Char} (1), we can only take $u = 0$ (if $u=1$, then $a = 2q^2>2^2q$) in the statement. But then we need to have $2<q<2^2$, so $q = 3$.   
\end{enumerate}

\section{Classification of distributed spectral pair in some cyclic groups}\label{classification}
We now consider  the distributed spectral pair of $\{0,1, \dots, M-1\}$ in the following cases: $M = p^k$ is a prime power, $M = pq$, and $M=p^2q$. We show that, in the first two cases, only spectral pairs are distributed spectral pairs; but in the last case, there are non-spectral but distributed spectral pairs. Throughout this section, we define for a subset $\A$ of non-negative integers, 
$$
P_{\mathcal A}(z) = \sum_{a\in{\mathcal A}}z^a
$$
$\Phi_n(x)$ denotes the cyclotomic polynomial of degree $n$, which is the minimal polynomial for the $n^{th}$ root of unity. We recall that the degree of $\Phi_n$, denoted by deg $\Phi_n(x)$, is equal to the Euler totient function $\phi(n) = n \prod_{p|n,\text{prime}}(1-\frac1p)$.
\begin{theorem}
    Suppose that $({\mathcal A}, {\mathcal B})$ is a distributed spectral pair in $\{0,1, \dots, M-1\}$ where $M = p^k$ is a prime power. Then $({\mathcal A}, {\mathcal B})$ is  a spectral pair in ${\mathbb Z}_M$. 
\end{theorem}

\begin{proof}
    Given the pair $({\mathcal A}, {\mathcal B})$ that is a distributed spectral pair, but not a spectral pair.  There exists $b_1,b_1'\in{\mathcal B}$ such that 
    $$
    \widehat{\nu}_{\mathcal A}\left(\frac{b_1-b_1'}{M^2}+ \frac{b_2-b_2'}{M}\right)  = 0 \ ~~\forall~ b_2, b_2'\in{\mathcal B}.
    $$
    We can write $b_1-b_1' = p^j r$ where $0\le j\le k-1$ and $r\in{\mathbb N}$ is relatively prime to $p$. Then
    $$
     \widehat{\nu}_{\mathcal A}\left(\frac{p^j r+ p^k(b_2-b_2') }{p^{2k}}\right) = 0 \ ~~\forall b_2, b_2'\in{\mathcal B}.
    $$ 
    We will take $b_2 = b_2'$. 
    This implies that 
    $$
    P_{\mathcal A}(e^{2\pi i r/p^{2k-j}}) =0
    $$
    showing that the cyclotomic polynomial $\Phi_{p^{2k-j}}(z)$ divides $P_{\mathcal A}(z)$. However, 
    $$
    \mbox{deg}(\Phi_{p^{2k-j}}) = (p-1)p^{2k-j-1}> p^k > \mbox{deg} (P_{\mathcal A}). 
    $$
    The factorization is impossible. Hence, the second equality never holds. The proof is complete.  
\end{proof}

The proof for two distinct primes will require the following theorem from De Bruijn. 

\begin{theorem}\label{de-B} (De Bruijn \cite{Bruijn1953})
 Let $n = p^{\lambda}q^{\mu}$ where $p,q$ are prime numbers and $P(x)$ is a polynomial with non-negative integer coefficients and deg $P(x)\le n$. Suppose $\Phi_n(x)$ divides $P(x)$.   Then there exist polynomials $A(x)$ and $B(x)$ with non-negative coefficients such that 
 $$
 P(x) =  A(x)\Phi_{p^{\lambda}} (x^{q^{\mu}})+ B(x)\Phi_{q^{\mu}} (x^{p^{\lambda}}).
 $$
\end{theorem}

We will call Theorem \ref{de-B} the De-Bruijn's theorem. We now prove the case for the product of two distinct primes.

\begin{theorem}
    Suppose that $({\mathcal A}, {\mathcal B})$ is a distributed spectral pair in $\{0,1, \dots, M-1\}$ where  $M = pq$ is a  product of two distinct primes with $p<q$. Then $({\mathcal A}, {\mathcal B})$ is  a spectral pair in ${\mathbb Z}_M$. 
\end{theorem}

\begin{proof}
Suppose that there exist $b,b'$ such that $P_{\mathcal A}(e^{2\pi i (b-b')/p^2q^2}) = 0$. Without loss of generality, we assume that $b>b'$.  We can write $b-b' = p^rq^sK$, where $K$ is relatively prime to $pq$. Note that $b-b'\in\{0,1,\cdots, pq-1\}$. In order not to exceed $pq-1$, it only happens that $b-b' = K, p^rK$ where $r$ is the largest integer such that $p^r<q$ or $qK$. In each case, we know that 
$$
\Phi_{p^2q^2}(x)|P_{\mathcal A}(x) \ \mbox{if $b-b' = K$},
$$  
$$
\Phi_{pq^2}(x)|P_{\mathcal A}(x) \ \mbox{if $b-b' = pK$},
$$
$$
\Phi_{q^2}(x)|P_{\mathcal A}(x) \ \mbox{if $b-b' = p^rK, r\ge 2$},
$$
$$
\Phi_{p^2q}(x)|P_{\mathcal A}(x) \ \mbox{if $b-b' = qK$}. 
$$
However, the first three cases are not possible because all cyclotomic polynomials have degrees larger than $pq$ and deg$(P_{\mathcal A})$ is at most $pq-1$. Indeed, recall that deg$(\Phi_n) = \phi(n)$, the Euler totient function, which has the formula
$$
\phi(n) = n \prod_{i=1}^k\left(1-\frac1{p_i}\right) \ \mbox{if} \ n = p_1^{\alpha_1}\cdots p_k^{\alpha_k}.
$$
Therefore, 
$$
\mbox{deg}(\Phi_{p^2q^2}) = pq (p-1)(q-1)>pq.
$$
The same holds true for $\Phi_{pq^2}(x)$ and $\Phi_{q^2}(x)$. Finally, in the last case, 
$$
\mbox{deg}(\Phi_{p^2q}) = p (p-1)(q-1).
$$
The expression $p (p-1)(q-1)<pq$ if and only if $p=2$. However, by the de-Bruijn's theorem, 
$$
P_{\mathcal A}(x) (\mbox{mod} \ x^{2^2q}-1) = \Phi_{2^2}(x^q) A(x)+\Phi_{q}(x^{2^2}) B(x)
$$
where $A,B$ are non-negative polynomials. As ${\mathcal A}\subset \{0,1,\cdots 2q-1\}$, $P_{\mathcal A}(x) (\mbox{mod} \ x^{2^2q}-1) = P_{\mathcal A}(x)$. We claim that $A(x), B(x)$ must be zero. Otherwise, if $x^a$ ($a>0$) exists in $A(x)$ or $B(x)$, then ${\mathcal A}$ will  contain elements $a+\{0,2q\}$ or $a+ \{0, 2^2\cdots (q-1)2^2\}$. In both cases, the largest element of  ${\mathcal A}$ is larger than $2q-1$, a contradiction. Hence, ${\mathcal A}$ is a singleton. This is not possible. The proof is complete.  
\end{proof}

We now proceed to the proof for $M = p^2q$. The proof strategy is similar to $M = pq$, but it has more complicated subcases. We first need a lemma.

\begin{lemma}\label{lemma4.4}
    Let $\A\subset \{0,1,\cdots,p^2q-1\}$ and $p<q$ are primes. Suppose that $\Phi_{pq^2}(x)$ divides $P_{\A}(x)$. Then $p=2$ and $q=3$ and $$
    \A = \{0,9\}, ~\{0,1,9,10\}, ~\{0,2,9,11\},\ \mbox{or} \ \{0,1,2,9,10,11\}.
    $$
\end{lemma}

\begin{proof}
Applying De-Bruijn's theorem, 
    \begin{equation}\label{eq:2a}
    P_{\mathcal A}(x) (\mbox{mod} \ x^{p^{}q^2}-1) = \Phi_{p^{}}(x^{q^2}) N(x)+\Phi_{q^2}(x^{p^{}}) M(x).
\end{equation}
As $p<q$, $pq(q-1)\ge p^2q$ and we must have $M(x)=0$. But  we need deg$\Phi_{p}(x^{q^2})=q^2(p-1)<p^2q$ which implies that $q(p-1)<p^2$. Since $p^2\le (p-1)(p+2)$, we have
$$
p<q<\frac{p^2}{p-1}\le p+2.
$$
Then $q=p+1$. Since $p$ and $q$ are primes, we must have $p=2$ and $q=3$. Now,  (\ref{eq:2a}) becomes
\begin{equation}\label{eq:2.1a}
    P_{\mathcal A}(x)  = \Phi_{2}(x^9) N(x) = (1+x^9)N(x).
\end{equation}
But $\A\subset \{0,1,\cdots, 11\}$. This forces that $N(x)$ has degree at most $2$.  This implies that $N(x) = 1, 1+x, 1+x^2$ or $1+x+x^2$. This completes the proof of the lemma. 
\end{proof}

\begin{theorem}
    Suppose that $({\mathcal A}, {\mathcal B})$ is a distributed spectral pair in $\{0,1, \dots, M-1\}$ where $M = p^2q$ and  $p,q$ are primes. Then it has to be one of the following cases:
    \begin{enumerate}
        \item $({\mathcal A}, {\mathcal B})$ is a spectral pair in ${\mathbb Z}_M$.
        \item $p=2, q=3$ and $({\mathcal A}, {\mathcal B})$ is a spectral pair in ${\mathbb Z}_{M^2}$. Moreover, ${\mathcal A}=\{0, 8\}$ and ${\mathcal B}=\{0,9\}$ after translations.
    \end{enumerate}
     
\end{theorem}
\begin{proof}
    Suppose that we do not have (1), we are going to study case by case that (2) is the only possible distributed spectral pair. As (1) fails,   there exist $b,b'$ (with $b>b'$) such that $P_{\mathcal A}(e^{2\pi i (b-b')/p^4q^2}) = 0$. We can write $b-b' = p^rq^sK$, where $K$ is relatively prime to $pq$. Note that $b-b'\in\{0,1,\cdots, p^2q-1\}$. It only happens that 
    $$b-b' = \begin{cases} 
    K; \\
    p^rK, \text{ where } p^{r-2}<q,;\\
    q^sK, \text{ where } q^{s-1}<p^2;\\
    pq^2K, \text{ where } q<p; \\
    pqK.
    \end{cases}$$
Note that the condition $p^{r-2}<q$ follows from the fact that we need to require $p^r<p^2q$ since $b-b'<p^2q$ (the same for the other cases). In each case, we know that 
$$
\Phi_{p^4q^2}(x)|P_{\mathcal A}(x) \ \mbox{if $b-b' = K$},
$$  
$$
\Phi_{p^{4-r}q^2}(x)|P_{\mathcal A}(x) \ \mbox{if $b-b' = p^rK$, $1\le r\le 3$},
$$
$$
\Phi_{q^2}(x)|P_{\mathcal A}(x) \ \mbox{if $b-b' = p^rK$, $r\ge 4$},
$$
$$
\Phi_{p^4q}(x)|P_{\mathcal A}(x) \ \mbox{if $b-b' = qK$},
$$
$$
\Phi_{p^4}(x)|P_{\mathcal A}(x) \ \mbox{if $b-b' = q^sK$, $s\ge 2$},
$$
$$
\Phi_{p^3}(x)|P_{\mathcal A}(x) \ \mbox{if $b-b' = pq^2K$, $q<p$},
$$
$$
\Phi_{p^3q}(x)|P_{\mathcal A}(x) \ \mbox{if $b-b' = pqK$}. 
$$
Since deg $\Phi_{p^4q^2}(x)$, deg $\Phi_{q^2}(x)$, deg $\Phi_{p^4q}(x)$ and deg $\Phi_{p^3}(x)$  are larger than or equal to $p^2q$ which is larger than deg $P_{\mathcal A}(x)$ (the last case because $q<p$), it remains to consider the following cases:  
$$b-b' = \begin{cases} 
     pqK, \text{ where } p=2;\\
     p^rK, \text{ where } p^{r-2}<q, 1\le r\le 3;\\
    q^2K, \text{ where } p(p-1)<q<p^2.\\
    \end{cases}$$
Indeed, the first two cases come from  $\Phi_{p^3q}(x)|P_{\A}(x)$ and $\Phi_{p^{4-r}q^2}(x)|P_{\A}(x)$ where $r\le 3$. In the last case, it follows from the fact that   $\Phi_{p^4}(x)|P_{\mathcal A}(x)$ and deg $P_{\A}(x)<p^2q$, so deg$\Phi_{p^4} = p^3(p-1)<p^2q$. This implies that $p(p-1)<q$. At the same time, $b-b' = q^s K<p^2q$ for some $s\ge 2$, so taking $K=1$ gives $q^{s-1}<p^2$.  Hence, $q \le q^{s-1}< p^2$ since $s\ge 2$.

In the following, we consider each case and deduce the complete structure of the distributed spectral pair. 

\medskip

\text{\bf Case 1}: $b-b' =2qK.$ \\

By the de-Bruijn theorem, we have 
\begin{equation}\label{eq:1}
    P_{\mathcal A}(x) (\mbox{mod} \ x^{2^3q}-1) = \Phi_{2^3}(x^q) N(x)+\Phi_{q}(x^{2^3}) M(x),
\end{equation}
where $N, M$ are non-negative polynomials. Clearly deg $P_{\mathcal A}(x) < 2^2q$, deg $\Phi_{2^3}(x^q)=2^2q$ and deg $\Phi_{q}(x^{2^3})=2^3(q-1)$ which implies that both deg $\Phi_{2^3}(x^q)$ and deg $\Phi_{q}(x^{2^3})$ are greater than deg $P_{\mathcal A}(x)$. It means $M(x)=N(x)=0$. It is a contradiction to (\ref{eq:1}). \\

\text{\bf Case 2}: $b-b' =p^rK, \text{ where } p^{r-2}<q, 1\le r\le 3.$ \\

We first claim that there are only three possible choices to $\A$, namely, 
\begin{equation}\label{eq_A}
\A = \{0,9\}, ~\{0,1,9,10\},~\{0,2,9,11\}, ~\{0,1,2,9,10,11\}.
\end{equation}
 To justify this claim, by the de-Bruijn theorem, we have 
\begin{equation}\label{eq:2}
    P_{\mathcal A}(x) (\mbox{mod} \ x^{p^{4-r}q^2}-1) = \Phi_{p^{4-r}}(x^{q^2}) N(x)+\Phi_{q^2}(x^{p^{4-r}}) M(x),
\end{equation}
where $N, M$ are non-negative polynomials. Note that 
\[
\mbox{deg}~ P_{\mathcal A}(x) < p^{2}q,~ \mbox{deg} ~\Phi_{p^{4-r}}(x^{q^2})=q^2p^{3-r}(p-1)~ \mbox{and deg} \ \Phi_{q^2}(x^{p^{4-r}})=p^{4-r}q(q-1).
\] 
Since $p^{4-r}q(q-1)\ge p^2q$ (as $p<q$), we must have $M(x)=0$. Therefore, $N(x)\not=0$. But for $r=1,2$, we have deg $\Phi_{p^{4-r}}(x^{q^2})=q^2p^{3-r}(p-1)\ge p^2q$. Thus, we can only have $r=3$ and  thus $\Phi_{pq^2}(x)$ divides $P_{\A}(x)$.  As we have $p<q$ in this case,   Lemma \ref{lemma4.4} shows that $p=2$ and $q=3$ and it also justifies  (\ref{eq_A}).

\medskip

We now analyze the structure of $\B$ for each case. Indeed, by (\ref{eq-dual}) and Theorem \ref{maintheorem11},  $(\A,\B)$ is a distributed spectral pair if and only if $(\B,\A)$ is a distributed spectral pair. Hence, for all $a\ne a'\in \A$, one of the following holds. 
\begin{itemize}
    \item [(i)] $P_{\mathcal B}(e^{2\pi i (a-a')/p^4q^2}) = 0$,
    \item [(ii)] $P_{\mathcal B}(e^{2\pi i (a-a')/p^2q}) = 0$.
\end{itemize}

\smallskip

\noindent\underline{\bf Subcase 1.  $\A = \{0,9\}$} We notice that $\A-\A = \{0,\pm9\}$ and if we have case (i), then $\B = \{0,8\}$. $(\A,\B)$ will be exactly the case (2) in the theorem. If we have Case (ii), we will have $\B = \{0,2\}$ and $(\A,\B)$ is a spectral pair, which is a contradiction, since we assume that we don't have a spectral pair in the beginning. 
\medskip

\noindent\underline{\bf Subcase 2.  $\A = \{0,1,9,10\}$ (or $\A = \{0,2,9,11\}$)} Note that $\A-\A = \{0,\pm1,\pm 8,\pm 9\}$ (or $\{0,\pm2, \pm 7,\pm 9,\pm 11\})$ and $\#\B = 4$. When $a-a' = 1$ (or 11), then (i) cannot hold since $\Phi_{p^4q^2}$ has degree larger than $p^2q$. so only $(ii)$ can hold and $\Phi_{p^2q}(x)$ divides $P_{\B}(x)$. Applying the De-Bruijn theorem, 
\begin{equation}\label{eq-P_B}
P_{\B}(x) = \Phi_{p^2}(x^q) M(x)+\Phi_q(x^{p^2}) N(x)
\end{equation}
where $M,N$ are non-negative polynomials with integer coefficients. Plugging in $x=0$, we have $4 = 2 M(0)+ 3 N(0)$. Hence, we must have $N(0) = 0$ and $M(0) = 2$. 
$$
P_{\B}(x) = \Phi_{2^2}(x^3) M(x) = (1+ x^6)M(x). 
$$
However, the largest element in $\B$ is at most 11. This forces $M$ is at most degree 5. Thus $M(x) = 1+x^j$ for some $1\le j\le 5$ and
$$
P_{\B}(x) = (1+x^6)(1+x^j).
$$

When $\A = \{0,1,9,10\}$, both (i) and (ii) cannot hold for $a-a' = 8$, which requires $\Phi_{2^33^2}(x)$ or $\Phi_{6}(x)$ divides $P_{\B}$. Therefore, this subcase does not contribute any distributed spectral pair. 

When $\A = \{0,2,9,11\}$, $P_{\A}(x) = \Phi_4(x)\Phi_2(x^9) = \Phi_4(x)\Phi_2(x)\Phi_6(x)\Phi_{18}(x)$. Now, $6\in \B$. Thus, we must require $e^{2\pi i 6/12}$ or $e^{2\pi i 6/(12)^2}$ to be a zero of $P_{\A}$. However, this means that $\Phi_2$ or $\Phi_{24}$ divides $P_{\A}$, which is not possible. 

\medskip

\noindent\underline{\bf Subcase 3.  $\A = \{0,1,2,9,10,11\}$} In this case, $\A-\A = \{0,\pm1,\pm2,\pm 7,\pm8,\pm9,\pm10,\\\pm11\}$ and hence it contains already all elements in subcase 1. The same argument in subcase 1 shows that $P_\B$ must satisfy (\ref{eq-P_B}) and $6 = 2M(0)+3N(0)$.  This implies that $M(0) = 0$ or $N(0) = 0$. This shows that 
$$
P_{\B}(x) = (1+x^6) M(x) \ \mbox{or} \ P_{\B}(x) = (1+x^4+x^8) N(x)
$$
However,  $P_{\B}(x) = (1+x^6) M(x)$ is impossible since deg$P_{\B}\le 11$, which means that $M(x)$ has at most degree $5$. $M(x)$ could have the form $1 + x^{m_1}+x^{m_2}$ where $1 \le m_1 < m_2 \le 5$. Using $a - a' = 2$, we
get that (ii) holds, obtaining $M(x) = 1 + x
^2 + x^4$ which leads to
$\B = \{0, 2, 4, 6, 8, 10\}$. Hence, $(\A,\B)$ must be a spectral pair. 

  If $P_{\B}(x) = (1+x^4+x^8) N(x)$, then $N(x) = 1+x$ or $1+x^2$ or $1+x^3$. We consider $a-a' = 9$. As (i) or (ii) holds, we have either $\Phi_{2^4}(x) = 1+x^8$ (if (i) happens) or $\Phi_{2^2}(x) = 1+x^2$ (if (ii) happens) divides $P_\B(x)$. The first is impossible since this will mean $P_{\B}$ has degree more than 12. Hence, we conclude that 
$\B= \{0,2,4,6,8,10\}$. In particular,  $(\A,\B)$ must be a spectral pair. Therefore, this subcase only contributes a spectral pair, which is not under our assumption in the beginning. 

\medskip

\text{\bf Case 3}: $b-b' =q^2K, \text{ where } p(p-1)<q<p^2.$ \\

The assumption that $P_{\mathcal A}(e^{2\pi i (b-b')/p^4q^2}) = 0$ implies that we have 
\begin{equation}\label{eq:3}
    P_{\mathcal A}(x)  = \Phi_{p^{4}}(x) N(x),
\end{equation}
where $N$ is a non-negative polynomial. By the fact that $(\B,\A)$ is also a distributed spectral pair, for all $a\ne  a'\in \A$, $P_{\mathcal B}(e^{2\pi i (a-a')/p^4q^2}) = 0$ or $P_{\B}(e^{2\pi i (a-a')/p^2q}) = 0$.

From (\ref{eq:3}), we know that there must exist $a,a'\in\A$ such that $a-a' = p^3$. If $P_{\mathcal B}(e^{2\pi i (a-a')/p^4q^2}) = 0$ for $a-a' =p^3$, we will have $\Phi_{pq^2}(x)$ divides $P_{\mathcal B}(x)$. By Lemma \ref{lemma4.4}, $p = 2$ and $q=3$ and $\B$ are of the form (\ref{eq_A}). This goes back to {\bf Case 2} and we know all the structures of $(\A,\B)$. 

Suppose now that $P_{\B}(e^{2\pi i (a-a')/p^2q}) = 0$ for $a-a' = p^3$. Then $\Phi_q(x)$ divides $P_{\B}(x)$. This implies that $q$ divides $\#\B = \#\A$. Hence, $pq$ divides $\#A$. Because we assume that $(\A,\B)$ is a non-spectral distributed spectral pair, there must be $a-a'$ such that $P_{\mathcal B}(e^{2\pi i (a-a')/p^4q^2}) = 0$. From all of our cases so far, we can only have $a-a' = q^2 K'$ for some integer $K'$ coprime to $p,q$. Hence, 
$$
P_{\B}(x) = \Phi_{q}(x)\Phi_{p^4}(x) M(x)
$$
for some non-negative integer polynomial $M$. Hence, $\B$ contains the following set:
$$
\{0,1,\cdots, q-1\}+ \{0,p^3,\cdots, p^3(p-1)\}. 
$$
Note that $p(p-1)<q<p^2$, 1 and $p$ are in $\B-\B$ and $p^3\in\B-\B$. It must satisfy the condition for the distributed spectral pair $(\A,\B)$. However, from the polynomial degree analysis and Lemma \ref{lemma4.4}, we must have $\Phi_{p^2q}, \Phi_{pq}, \Phi_q$ divides $P_{\A}$. Hence, 
$$
P_{\A}(x) = \Phi_{p^{4}}(x) \Phi_{q}(x^{p^2})N'(x),
$$
We now compute the 
$$
\mbox{deg} ~ \Phi_{p^{4}}(x) \Phi_{q}(x^{p^2}) = p^3(p-1)+p^2(q-1) = p^2q+ p^2(p^2-p-1)>p^2q.
$$ 
This is a contradiction. Hence, Case 3 does not result in any non-trivial distributed spectral pair. This completes the whole proof.
\end{proof}


\section{Discussions and Open problems}\label{discussion}

\subsection{Discrete settings.} The classification of distributed spectral pair in $\{0,1,\cdots M-1\}$ suggest the following conjecture may be correct. 
\begin{conj}\label{conj}
    Let $(\A, \B)$ be a distributed spectral pair in $\{0,1,\cdots M-1\}$. Then $(\A, \B)$ is either a spectral pair in $\Z_M$ or $\Z_{M^2}$.
\end{conj}
So far, we didn't find any counter-example of Conjecture \ref{conj}. Moreover, supposing Conjecture \ref{conj} holds, if $(\A, \B)$ is a distributed spectral pair in $\Z_M$ but not a spectral pair in $\{0,1,\cdots M-1\}$, then by Theorem \ref{theorem-spectral-m2}, the set $\A$ contains only two elements. 

{\it Added Note: During the revision of the article, we discovered that Conjecture \ref{conj} is not true. The example is as follows: $M = 24$, $A = \{0,4,16,20\}$ and $\B = \{0,3,18,21\}$. $(\A,\B)$ cannot be a spectral pair in $\Z_{24}$ and it cannot be a spectral pair in $\Z_{(24)^2}$ using Proposition \ref{prop_spectral_M2}. Yet, a direct calculation shows that this is a distributed spectral pair. }

On the other hand, if we drop the condition that both $\A$ are $\B$ in $\{0,1\cdots, M-1\}$, then we can actually find a pair $(\A, \B)$ satisfies the requirement (2) in Theorem \ref{maintheorem1} but it is neither a spectral pair in $\Z_M$ or $\Z_{M^2}$ and it has more than two elements. We would like to thank Professor Li-Xiang An for pointing this out to us.  The example is as follows:

\begin{example}
    Let $M = 2^3\cdot m$, $\A = \{0,1, 2^4, 2^4+1\}$ and $\B = \{0, b,2b,3b\}$, where $b = 2m^2$ and $m$ is an odd integer strictly larger than 1. Then $(M, \A,\B)$ satisfies the condition (2) in Theorem \ref{maintheorem1}. Notice that in this case,  $\B$ is not in $\{0,1,\cdots, M-1\}$.
\end{example}

\begin{proof}
    We first factorize the polynomial
    $$
    P_{\A}(x) = \Phi_2(x)\Phi_{2^5}(x).
    $$
    The difference set of $\B$ is equal to $
    \B-\B = \{0, \pm b, \pm2b,\pm 3b\}$.    A direct check shows that 
    $$
    \widehat{\nu_{\A}} \left(\frac{\pm2b}{M}\right) = 0,  ~   \widehat{\nu_{\A}} \left(\frac{\pm b}{M^2}\right)=    \widehat{\nu_{\A}} \left(\frac{\pm3b}{M^2}\right)  = 0.
    $$
    By a further check, we see that 
    $$
    \widehat{\nu_{\A}} \left(\frac{\pm j b}{M^2} + \frac{\pm  k b}{M}\right)=0
    $$
    for all $j = 1,3$ and $k = 0, 1,2,3$. This fulfills the definition of a distributed spectral pair. Finally,  $\widehat{\nu_{\A}} \left(\frac{\pm2b}{M^2}\right)\ne 0$ and $\widehat{\nu_{\A}} \left(\frac{\pm b}{M}\right)\ne 0$, so $(\A,\B)$ cannot be a spectral pair in $\Z_M$ and $\Z_{M^2}$. This completes the proof of this example. 
\end{proof}

Indeed, this example was first obtained in \cite{AL2023}. It can be generalized to more general alphabet sets called the product-form Hadamard triple. However, we did not find any such examples that  can lie completely inside $\{0,1,\cdots M-1\}$. 

When an alphabet set $\B$ lies outside $\{0,1,\cdots M-1\}$, $\B_k$ will lie outside of $\{0,1\cdots M^k-1\}$. If we let $\widehat{\B}_k  = \{b_k (\mbox{mod} \ M^k): b_k\in\B_k\}$, then $\widehat{\B}_k\subset \{0,1,\cdots, M^k-1\}$ and we see that $(M^k,\A_k,\widehat{B}_k)$ also attains the most uncertain exponent. While the self-similarity structure of $\widehat{B}_k$ is destroyed, we still obtain an uncertainty exponent using $\widehat{B}_k$.

\subsection{Continuous settings} Our study also shed some light on the most uncertain situation in the continuous setting. The {\it semi-classical Fourier transform} with parameter $h>0$,
$$
{\mathcal F}_h: L^2(\R)\to L^2(\R), ~~~ {\mathcal F}_hf(\xi)  = \frac{1}{\sqrt{2\pi  h}} \int f(x) e^{- i \frac{\xi \cdot x}{h}} d\xi.
$$
For a given compact subset $X\subset \R$, we let $X_h = X+[-h,h]$ to be the $h$-neighborhood of $X$.  Suppose that $X, Y$ are compact subsets of $\R$.  It is not hard to show that 
\begin{equation}\label{eq-trivial}
\|{\mathbf 1}_{X_h} {\mathcal F}_h{\bf 1}_{Y_h}\|_{L^2(\R)\to L^2(\R)} \le 1
\end{equation}
and 
\begin{equation}\label{eq-volume}
\|{\mathbf 1}_{X_h} {\mathcal F}_h{\bf 1}_{Y_h}\|_{L^2(\R)\to L^2(\R)} \le (2\pi h)^{-1/2} \sqrt{m(X_h)\cdot m(Y_h)}    
\end{equation}
where $m$ denotes the Lebesgue measure on $\R$. The inequality is due to the fact that they are unitary operators, while the second upper bound is obtained by computing its Hilbert-Schmidt norm. We say that $X,Y$ are {\bf $(\delta,C)$-Ahlfors–David regular}  if  there exists a measure $\mu$ supported on $X$ such that 
$$
C^{-1}r^{\delta }\le \mu(( x-r,x+r)) \le C r^{\delta}, ~ \forall x\in X.
$$
We know that $m(X_h) \asymp  h^{1-\delta}$ (since the box dimension of $X$ is $\delta$). Thus,  (\ref{eq-trivial}) and (\ref{eq-volume}) implies that 
$$
\|{\mathbf 1}_{X_h} {\mathcal F}_h{\bf 1}_{Y_h}\|_{L^2(\R)\to L^2(\R)}\lesssim h^{\max\{0,\frac12-\delta\}}.
$$ 
Bourgain, Dyatlov and Jin \cite{BourgainDyatlov2018, DJ2018} improved the trivial bound above as follows:

\begin{theorem}
There exists $\varepsilon_0>0$ depending only $\delta$ and $C>0$ such that 
\begin{equation}\label{eq_FUP}
 \|{\mathbf 1}_{X_h} {\mathcal F}_h{\bf 1}_{Y_h}\|_{L^2(\R)\to L^2(\R)}\lesssim h^{\max\{0,\frac12-\delta\}+\varepsilon_0}   
\end{equation}
for all $(\delta,C)$-Ahlfors–David regular sets $X,Y$. 
\end{theorem}
We will say that {\bf fractal uncertainty principle holds for $(X_h,Y_h)$ with exponent $\beta$} if (\ref{eq_FUP}) holds with $h^{\beta}$ on the right. Note that the continuous version immediately implies the discrete version in Theorem \ref{th_DJ}.  Fractal uncertainty principles also exist for operators with general non-linear phase functions and it was also recently generalized into high dimensions \cite{C2024}.

It is possible to obtain a lower bound for the continuous FUP setting. Notice that the lower bound would not hold if we just consider $X,Y$ are $\delta$-regular sets up to scale $h$ in \cite{BourgainDyatlov2018}, because $X,Y$ are now allowed to be a measure zero set, so ${\bf 1}_{X} {\mathcal F}_h{\bf 1}_{Y} f =0$ for all $f\in L^2(\R)$. We may however consider the neighborhood of $\delta$-regular sets. The following proposition shows that $\frac{1-\delta}{2}$ is still the most uncertain exponent in the continuous setting. 

\begin{proposition}
Let $X,Y$ be $\delta$-Ahlfors regular sets and $X_h,Y_h$ be its $h$-neighborhood. Then 
$$
\|{\bf 1}_{X_h} {\mathcal F}_h{\bf 1}_{Y_h}\|_{L^2(\R)\to L^2(\R)} \gtrsim h^{\frac{1-\delta}{2}}
$$
where the implicit constant depends on the set $X$. 
\end{proposition}

\begin{proof}
Since $Y_h$ contains $(y-h,y+h)$ for some $y\in Y$. We can now take $f = {\bf 1}_{(y-h,y+h)}$. Then 
$$
{\bf 1}_{X_h} {\mathcal F}_h{\bf 1}_{Y_h} f(\xi) = {\bf 1}_{X_h}(\xi)\frac{1}{\sqrt{2\pi h }} \cdot \int_{y-h}^{y+h} e^{-i\frac{\xi x}{h}} ~ dx  = {\bf 1}_{X_h}(\xi)\cdot {\sqrt{\frac{h}{2\pi} }} \cdot e^{i \xi \frac{y}{h}}\cdot \int_{-1}^{1 } e^{-i\xi x} ~ dx  
$$
by a change of variable $x\to \frac{x-y}{h}$. As $\int_{-1}^{1 } e^{-i\xi x} ~ dx  = \frac{2\sin \xi}{\xi} = 2\mbox{sinc}(\xi)$, we have 
$$
\|{\bf 1}_{X_h} {\mathcal F}_h{\bf 1}_{Y_h}\|_{L^2\to L^2}\ge \frac{\|{\bf 1}_{X_h} {\mathcal F}_h{\bf 1}_{Y_h} f\|_{L^2(\R)}}{\|f\|_{L^2(\R)}}  = \frac{1}{\sqrt{\pi}}\cdot \left(\int_{X_h} \left |\mbox{sinc} (\xi)\right|^2d\xi\right)^{1/2}.
$$
But $\mbox{sinc}(\xi)$ is a continuous function whose zeros are precisely at the non-zero integer multiples of $\pi$. As $X$ is $\delta$-Ahlfors regular, we can always find a point $t\in X\setminus \pi \Z$, $r>0$ and $X'=X\cap (t-r,t+r)$ so that $X'$ does not contain points in $\pi \Z$ and $X'$ is also $\delta$-Ahlfors regular. Hence,  for some constant $c>0$, depending only on $t,r$, 
$$
\left(\int_{X_h} \left |\mbox{sinc} (\xi)\right|^2d\xi\right)^{1/2}\ge c \cdot (m(X'_h))^{1/2}.
$$
Since $X'$ is $\delta$-Ahlfors regular, $m(X'_h)\gtrsim h^{1-\delta}$. This shows that our desired conclusion holds. 
\end{proof}

It now makes sense to discuss what pair $(X,Y)$ of $\delta$-regular sets attain the most uncertain exponent. A measure is called a {\bf spectral measure} if there exists a countable set $\Lambda \subset \R$ such that $\{e^{2\pi i \lambda x}\}_{\lambda\in \Lambda}$ forms an orthonormal basis of $L^2(\mu)$. The first singular spectral measure was discovered by Jorgensen and Pedersen \cite{JP98}, who show that the middle-fourth Cantor measure is a spectral measure. Indeed, {\L}aba and Wang \cite{LW02} showed that   $(\A,\B)$ is a spectral pair, the associated self-similar measure $\mu(M,\A)$ and $\mu(M,\B)$ (c.f. \ref{eq-self-similar}) are all spectral measures. For example, if $M = 4$, $\A = \{0,2\}$ and $\B = \{0,1\}$, then $\mu (4,\{0,2\})$ is a spectral measure with 
$$
\Lambda = \bigcup_{k=1}^{\infty} \B_k = \left\{\sum_{j=0}^{k-1}4^j\varepsilon_j: \varepsilon_j\in\{0,1\},~ k\in\N\right\}.
$$
{\L}aba-Wang's result was generalized to higher dimensions by Dutkay, Haussermann, and the first-named author \cite{DHL2019} (see \cite{DLW2017} for a survey about these results). A type of product-form spectral pair was obtained in \cite{AL2023} and they can be used to produce distributed spectral pairs. 

\medskip

For a distributed spectral pair $(\A, \B)$ in $\Z_M$, we call $\mu(M,\A)$ and $\mu (M,\B)$ the {\bf distributed spectral measure}. From our results on the discrete Cantor sets, it seems plausible to ask the following question:
\begin{ques}\label{question}
    Does the support of a distributed spectral measure $\mu(M,\A)$ and $\mu(M,\B)$ with $(\A,\B)$ are distributed spectral pair in $\Z_M$  attain the most uncertain exponent?
\end{ques}
We remark that the affirmative answer of Question \ref{question} will be a continuous version of Theorem \ref{maintheorem1}. So far, we don't know the answer of Question \ref{question}. But we believe that the answer would be yes.

\subsection*{Acknowledgment} The authors would like to thank Professors Mihalis Kolountzakis  and Li-Xiang An for some useful discussions. Part of the work was done when the authors were visiting the Chinese University of Hong Kong (CUHK) and Central China Normal University (CCNU). They would like to thank Professor De-Jun Feng, CUHK Department of Mathematics, Li-Xiang An and CCNU for their hospitality. We would like to thank the referee for a careful reading of the paper and providing many constructive comments. 

Chun-Kit Lai is partially supported by the AMS-Simons Research Enhancement Grants for Primarily Undergraduate Institution (PUI) Faculty.  Ruxi Shi was supported by NSFC No. 12571198 and No. 12231013. This work has also been supported by the New Cornerstone Science Foundation through the New Cornerstone Investigator Program.

\bibliographystyle{alpha}
\bibliography{universal_bib}

\end{document}